\documentclass[12pt]{amsart}
\usepackage{amsmath, amsthm, latexsym, amssymb, amsfonts}
\usepackage[pdftex]{graphicx,color} 

\newenvironment{pf}{\proof[\proofname]}{\endproof}
\theoremstyle{plain}
\newtheorem{Th}{Theorem}[section]
\newtheorem{Cor}[Th]{Corollary}
\newtheorem{Prop}[Th]{Proposition}
\newtheorem{Lemma}[Th]{Lemma}
\numberwithin{equation}{section} \theoremstyle{definition}

\newtheorem{Ex}{Example}

\newtheorem{Def}[Th]{Definition}

\usepackage{mathtools}
\DeclarePairedDelimiter\floor{\lfloor}{\rfloor}

\newcommand{\cal}[1]{\mathcal{#1}}

\newcommand{\R}{\mathbb R}

\newcommand{\cD}{{\mathcal D}}
\newcommand{\cB}{\cal B}

\newcommand{\tdet}{{\rm tdet\,}}
\newcommand{\tropdet}{{\rm tropdet\,}}

\newcommand{\rl}[1]{Lemma~\ref{L:#1}}
\newcommand{\rp}[1]{Proposition~\ref{P:#1}}

\newcommand{\rex}[1]{Example~\ref{ex:#1}}
\newcommand{\re}[1]{(\ref{e:#1})}
\newcommand{\rc}[1]{Corollary~\ref{C:#1}}
\newcommand{\rt}[1] {Theorem~\ref{T:#1}}
\newcommand{\rd}[1]{Definition~\ref{D:#1}}
\newcommand{\rf}[1]{Figure~\ref{F:#1}}

\begin{document}
\author{Sailaja Gajula}
\address{Department of Mathematics, Indiana University, Bloomington, IN 47405}
\email{sgajula@indiana.edu}
\author{Ivan Soprunov} 
\address{Mathematics Department, Cleveland State University, Cleveland, OH 44115}
\email{i.soprunov@csuohio.edu}
\author{Jenya Soprunova}
\address{Department of Mathematical Sciences,
Kent State University, Kent, OH 44242}
\email{soprunova@math.kent.edu}
%\date{September 23, 2003}

% \thanks{The second author is partially supported by NSA Grant H98230-13-1-0279}

\title{Tropical Determinant on Transportation Polytopes}

\begin{abstract}  
Let ${\cD^{k,l} }(m,n)$ be the set of all the integer points in the transportation polytope of 
$kn\times ln$ matrices with row sums $lm$ and column sums $km$.  
In this paper we find the sharp lower bound on the tropical determinant over the set ${\cD^{k,l}}(m,n)$. 
This integer piecewise-linear programming problem in arbitrary dimension turns out to be equivalent  to an integer non-linear (in fact, quadratic) optimization problem in dimension two.
We also compute the sharp upper bound on a modification of the tropical determinant, where 
the maximum over all the transversals in a matrix is replaced with the minimum.
\end{abstract}
\maketitle
%%%%%%%%%%%%%%%%%%%%%%%%%%%%%%%

\section*{Introduction}
In this paper we generalize the results of  ``Tropical determinant of integer doubly stochastic matrices'' 
\cite{DHS}  to the class of all rectangular integer matrices with fixed row and column sums. The discussion in \cite{DHS} started with cheater's Rubik's cube problem: When solving  Rubik's cube by peeling off and replacing stickers, how many stickers do we need to peel off and replace in the worst case scenario? This problem generalizes to a very natural sorting question: Assume that we have $n$ pails with $m$ balls in each.  Each ball is colored in one of $n$ colors and we have $m$ balls of each color.  What is the smallest  number of balls  we need to move from one pail to another in the worst case scenario so that the balls are sorted by color?

This problem turns out to be equivalent to finding the sharp lower bound on the tropical determinant of integer 
matrices $A=(a_{ij})$ of given size $n$ with given row and column sums $m$. 
To see this, let the entry $a_{ij}$
be equal to the number of  balls  of color $i$ in pail  $j$.  We would like to assign each pail a color so that the overall number of balls that we need to move is the smallest possible. That is, we would like to find a transversal of $A$ with 
the largest possible sum of entries, which is the definition of the tropical determinant $\tdet A$ of $A$. 

The set of all (real) doubly stochastic  $n\times n$ matrices forms a convex polytope  in $\R^{n^2}$,
the Birkhoff polytope $\cB_n$ (see \cite{Birk}). The set of integer $n\times n$ matrices with 
row and column sums equal $m$ can then be identified with the set of integer points of its
$m$-dilate $m\cB_n$. The tropical determinant is a piece-wise linear function on $m\cB_n$. Therefore, the described problem is equivalent to minimizing  this function over the integer points of the polytope, i.e. solving an integer piecewise-linear programming problem.  This was done in \cite{DHS}.

In the current paper we are working on a natural generalization of this problem, where we replace the Birkhoff polytope with any transportation polytope. A transportation polytope is a convex polytope consisting of nonnegative rectangular matrices of given size with fixed row and column sums. The set of integer such matrices 
is identified with the set of integer points of a transportation polytope.
Our goal is to compute the sharp lower bound for the tropical determinant on 
integer points of a transportation polytope. Surprisingly, this integer piecewise-linear programming problem in arbitrary dimension reduces to an integer non-linear (in fact, quadratic) optimization problem in dimension two (see \rt{mainshort}).

This problem has a similar combinatorial interpretation. Suppose there are $R$ balls of each of $t$ different colors, totaling $tR$ balls. Suppose they are placed into $s\geq t$ different pails with $C$ balls in each pail (so
$sC=tR$). We want to sort the balls by color in some $t$ of the $s$ pails, by replacing balls from one pail to another. 
What is the smallest  number of balls  we need to move from one pail to another to achieve this
in the worst case scenario? Similar to above,
let $a_{ij}$ be the number of  balls  of color $i$ in pail  $j$. We obtain an $r\times s$ matrix $A=(a_{ij})$ 
whose row sums are $R$ and column sums are $C$. The smallest number of moves to sort the balls 
is then $tR-\tdet A$. Thus, to answer the above question one needs to find the  sharp lower bound for the tropical determinant over all such matrices $A$.

In this paper we build on the methods developed in \cite{DHS}. We were able to simplify the arguments to the point where the desired generalization became possible. Also the answer in the general setting is more transparent.
Our methods are elementary and do not rely on other results except for Hall's Marriage theorem.

Following \cite{DHS} we also consider and solve a version of the problem where in the definition of the tropical determinant the minimum over all the transversals is replaced with the maximum. In this case we are 
interested in the sharp upper bound over the integer points of the transportation polytope.  As in \cite{DHS}, this version of the problem turns out to be significantly easier than the problem we start with. 

In 1926 van der Waerden conjectured that the smallest value of the permanent of $n\times n$ doubly stochastic (with row and column sums equal to one)
matrices is attained on the matrix all of whose entries are equal to $1/n$, and this minimum is attained only once. 
This conjecture was proved independently by 
 Egorychev \cite{Eg}  and Falikman \cite{Fal}  in 1979/80. In \cite{BuBu} Burkard and Butkovich proved a tropical version of the conjecture, where the permanent is replaced with the tropical determinant. Results of this paper and \cite{DHS} provide an integral tropical version of the van der Waerden conjecture.

%While in \cite{DHS}, the story develops from simple cases to harder ones, in this paper, we first prove our main result 
%(\rt{mainshort}) and then use it to understand simpler cases as well as recover
%results from \cite{DHS}. 

%%%%%%%%%%%%%%%%%%%%%%%%%%%%%%%%%%%%%%%%%%%%%%%%%%%%%%%%%%%%%%%%%%%%%%%%%%
\vspace{2cm}
\section{Definitions}
 Let $A=(a_{ij})$ be an $nk$ by $nl$ matrix where $\gcd(k,l)=1$ and $a_{ij}$ are non-negative integers. Let all the row sums in $A$ be equal to $a$ and all the column sums be equal to $b$.
Computing the sum of all the entries in $A$  in two different ways we get  $ka=lb$, which implies  $a=ml$, $b=mk$ for some integer $m$.

\begin{Def} Let $k\leq l$. Define  $\cD^{k,l}(m,n)$ to be the set of all $nk\times nl$ matrices with nonnegative integer entries whose row sums are $ml$ and columns sums are $mk$.
\end{Def}

\begin{Def} For an $s\times t$ matrix $A=(a_{ij})$ with $s\leq t$, its transversal $T$ is a set $\{a_{1i_1},\dots, a_{si_s}\}$  where $1\leq i_1<\dots< i_s\leq t$. 
Furthermore, let $|T|=a_{1i_1}+\cdots+a_{si_s}$ and $\mathcal{T}(A)$ be the set of all transversals of $A$.  For $t\geq s$ we define transversals of $A$ to be transversals of its transpose $A^{T}$. 
\end{Def}

\begin{Def} The {\it tropical determinant} of a matrix $A=(a_{ij})$ is
$$\tdet(A)=\max_{T\in\mathcal{T}(A)}|T|.$$
We will refer to a transversal of $A$ on which this maximum is attained as a {\it maximal transversal} of $A$.
\end{Def}
Clearly, the set of transversals and, hence, the tropical determinant are invariant under row and column swaps of $A$.

Let $L^{k,l}(m,n)$ denote the sharp lower bound on the tropical determinant over the set $\cD^{k,l}(m,n)$, that is,
$$L^{k,l}(m,n)=\min_{A\in \cD^{k,l}(m,n)} \tdet(A).
$$
Our main goal in this paper is to compute $L^{k,l}(m,n)$.

\begin{Ex}\label{ex:tdet9} Let $n=5$, $k=1$, $l=2$, and $m=6$. Then the matrix
\[A=\left(
\begin{array}{cccccccccc}
0&0&2&1&\fbox{2}&1&2&1&2&1\\
0&0&1&\fbox{2}&1&2&1&2&1&2\\
2&2&\fbox{1}&1&1&1&1&1&1&1\\
2&\fbox{2}&1&1&1&1&1&1&1&1\\
\fbox{2}&2&1&1&1&1&1&1&1&1\\
\end{array}
\right)
\]
lies in $\cD^{1,2}(6,5)$. The boxed elements form a maximal transversal of $A$. Thus $\tdet A=9$.
We will later show that $L^{1,2}(6,5)=9$, that is, the minimum of the tropical determinant on $\cD^{1,2}(6,5)$  is attained at this matrix. 
\end{Ex}

One of our tools is Hall's marriage theorem and, following \cite{DHS}, we restate this theorem and its simple corollaries here, making a small adjustment to the case of rectangular matrices.
The theorem in our formulation deals with a maximal zero submatrix of $A$, that is a zero submatrix of $A$ whose sum of dimensions is the largest possible.

\begin{Th}[Philip Hall \cite{Hall}]\label{T:Hall} Let $A$ be an $s\times s$ 0-1 matrix. Then there is a transversal in $A$ that consists of all 1's if and only if a maximal  zero submatrix in $A$ has sum of dimensions less than or equal to $s$. 
\end{Th}

For our future discussion we will need the following two corollaries.

\begin{Cor}\label{C:marriage} Let $A$ be an $s\times t$ 0-1 matrix. Then there is a transversal  in $A$ that consists of all 1's if and only if a maximal zero submatrix of  $A$  has sum of dimensions less than or equal to $\max(s,t)$. 
\end{Cor}

\begin{pf} Let us assume that $s\leq t$. Extend $A$ to a square 0-1 matrix by appending to $A$ $t-s$ rows consisting of all $1$'s and apply Hall's marriage theorem to the resulting matrix.
\end{pf}

Let $A$ be an $s\times t$ 0-1 matrix and $W$ be  a maximal  zero submatrix of  $A$. Then after some row and column swaps $A$ can be written in the form
\[
A=\left(\begin{array}{cc}   
X&Y\\
Z&W
\end{array}\right).
\]

\begin{Cor}\label{C:Marriage}  Both $Y$ and $Z$ have a transversal that consists of all 1's.
\end{Cor}
\begin{pf} 
Let $W$ be $d_1\times d_2$  and a maximal zero submatrix  of $Y$ be $s_1\times s_2$. We can assume that it is in the lower right corner of $Y$, right on top of  $W$. Then the lower right $(s_1+d_1)\times s_2$ block of $A$ consists of all zeroes which implies  $s_1+d_1+s_2\leq d_1+d_2$, and so $s_1+s_2\leq d_2$. By \rc{marriage}  there exists a  transversal in $Y$ that consists of all $1$'s. Similarly, such a transversal exists in $Z$.
\end{pf}

\section{Bound on $L^{k,l}(m,n)$}

We start with two simple observations concerning the tropical determinant of an arbitrary matrix.

\begin{Lemma}\label{L:easy} Let $B$ be an $s\times t$ matrix with $s\geq t$. Then $\tdet B$ is at least the sum of all the entries in $B$, divided by $s$. In particular, if all row sums of $B$ are bounded from below by $b$, then $\tdet B\geq b$.
\end{Lemma}
\begin{pf} The set of entries of $B$ can be partitioned into $s$ transversals $T_1,\dots,T_s$. 
Since $|T_i|\leq\tdet B$,  the sum of all  entries of $B$ does not exceed $s\,\tdet B$.
\end{pf}

\begin{Lemma}\label{L:key} Let $Q$ be an $s\times t$ matrix with $s\leq t$. Let $a$ be any element in a maximal transversal of $Q$.
Then $$R+C\leq \tdet Q +sa,$$ where $R$ and $C$ are the sum of entries in the row and column that contain~$a$.
\end{Lemma}
\begin{pf} After necessary row and column swaps we can assume that $a=a_s$ is in the position $(s,s)$ 
and that  $\tdet Q=a_1+\dots+a_s$, where
$$
 Q=\left( \begin{array}{cccccc}
a_1   & & & c_1     &  &   \\
  & a_2 &          &c_2         &  &   \\
  &        & \ddots     &\vdots  & &     \\
b_1&b_2 &\dots&  a_s &\dots   &b_t    \\
\end{array} \right),
$$
and $C=c_1+\cdots+c_{s-1}+a_s$ is the $s$-th column sum, while $R=b_1+\cdots+b_{s-1}+a_s+b_{s+1}+\cdots+b_t$ is the $s$-th row sum.

We have $b_j+c_j\leq a_j+a_s$ for $j=1,\dots, s-1$ since otherwise we could switch columns $j$ and $s$ in $Q$ to get a larger transversal. We also have $b_j\leq a_s$ for $j=s+1,\dots, t$. Summing these up over $j=1,\dots,t$ we get 
$$R+C\leq \tdet Q+sa_{s}.$$
\end{pf}

Recall that  $A=(a_{ij})\in\cD^{k,l}(m,n)$ is an $nk\times nl$ matrix where $k\leq l$, $\gcd(k,l)=1$,  and $a_{ij}$ are non-negative integers. The row sums of $A$ are equal to $ml$ and the column sums are equal to $mk$.

Now divide $m$ by $n$ with remainder, $m=qn+r$, where $0\leq r<n$. Let $W$ be a submatrix of $A$ with entries less than or equal to $q$ with the largest sum of dimensions.
Then after some column and row swaps $A$ can be written in the form
\[
A=\left(\begin{array}{cc}   
X&Y\\
Z&W
\end{array}\right).
\]
Let  $X$ be of size $t_1$ by $t_2$.

\begin{Lemma}\label{L:lknr} We have
$$qt_1t_2+r(t_1l+t_2k)\geq klnr.
$$
\end{Lemma}

\begin{pf}
Let $\Sigma_W$ and $\Sigma_Y$  be the sums of all the entries in blocks $W$ and  $Y$. Then  $\Sigma_W\leq q(nk-t_1)(nl-t_2)$ since all the entries of $W$ do not exceed $q$. Hence 
$$\Sigma_Y= (nl-t_2)mk-\Sigma_W\geq (nl-t_2)mk-q(nk-t_1)(nl-t_2).
$$
On the other hand,  $\Sigma_Y\leq t_1ml$. Putting these two inequalities together we get
$$(nl-t_2)mk-q(nk-t_1)(nl-t_2)\leq  t_1ml,
$$
which is easily seen to be  equivalent to $qt_1t_2+r(t_1l+t_2k)\geq klnr$ using $m=qn+r$.
This argument also shows that $qt_1t_2+r(t_1l+t_2k)-klnr$ is an upper bound for $\Sigma_X$. 
\end{pf}

This lemma motivates the following definition. 
\begin{Def}\label{D:xy} Let $x$ and $y$ be integers satisfying  $x\geq rk$, $y\geq rl$, and
\begin{equation}\label{e:xy}
qxy+r(xl+yk)\geq klnr,
\end{equation} 
whose sum $x+y$ is the smallest possible.
\end{Def}
Note that while $x+y$ is defined uniquely, this is not necessarily true for $x$ and~$y$.
Also, the conditions $x\geq rk$ and $y\geq rl$ will be necessary for the construction
in \rp{constr2}.

Recall that 
\[
A=\left(\begin{array}{cc}   
X&Y\\
Z&W
\end{array}\right),
\]
where $W$ is a maximal submatrix that consists of elements not exceeding $q$ and $X$ is of size $t_1$ by $t_2$.

\begin{Lemma}\label{L:main} Let $t_1+t_2\leq nk$. Then
$$
\tdet A\geq\min\left(nk(q+1), \tdet Y+\tdet Z+(nk-t_1-t_2)q\right).
$$
\end{Lemma}

\begin{pf}  
Consider the set of all transversals in $A$ which contain a maximal transversal in $Y$ and
a maximal transversal in $Z$. Choose one such transversal $T_A$ with the largest sum $|T_A|$.
Let  $T_Y\subset T_A$ and $T_Z\subset T_A$ be the corresponding maximal transversals in
$Y$ and $Z$, respectively. 
Note that since $t_1+t_2\leq nk$, the  transversals in $Y$ and $Z$ have respectively $t_1$ and $t_2$ entries.
Cross out the rows and columns of $A$ which contain 
$T_Y$ and $T_Z$ to get  an $(nk-t_1-t_2)\times(nl-t_1-t_2)$ submatrix $Q$ of $W$.
Then $T_A=T_Y\cup T_Z\cup T_Q$, where the transversal $T_Q$ of $Q$ is also maximal
(by construction). Therefore, we obtain
\begin{equation}\label{e:0}
\tdet A\geq |T_A|=\tdet Y+\tdet Z+\tdet Q.
\end{equation}
First, assume that $Q$ contains a transversal all of whose elements are equal to $q$. Then
we have $\tdet Q=(nk-t_1-t_2)q$ and the statement follows from the above inequality.

Next, assume that every maximal transversal of $Q$ has an entry less than or equal to $q-1$.
We can rearrange the rows and columns of $A$ as follows
$$
 A=\left( \begin{array}{cccc|cccccc|cccc}
&  &  & &        &      &  &c_1   & &         &a_1&&&\\
&  & X & &        &      &  &c_2    &  &      &     &a_2&&\\
&  &  & &        &      &  &\vdots      & &    &     &&\ddots&\\
&  &  & &        &      &   &c_{t_1}    &  &   &     &&& a_{t_1}\\
\hline
&  &  & &e_1   &                  &  &f_1    &  &       & &&&\\
&  &  & &        &e_2 &          &f_2         &  &   &     &&&\\
&  &  & &        &      &\ddots  &\vdots    &  &     &     &&&\\
d_1&d_2 &\dots  &d_{t_2} &g_1   &g_2 &\dots   &e_{t_3}  &\dots &g_{t_4}       &i_1     &i_2&\dots& i_{t_1}\\
\hline
b_1&      &          &            & & &  &h_1  & &          & &&&\\
     &b_2 &          &            & & &  &h_2      &  &         &     &&&\\
     &      &\ddots  & &        & & &\vdots              &&&\\
     &      &\          &b_{t_2}& &  &  &h_{t_2}    & &      &    &&&

\end{array} \right)
$$
Here the middle block is $Q$,
$T_Y=\{a_1,\dots, a_{t_1}\}$, $T_Z=\{b_1,\dots, b_{t_2}\}$, and $T_Q=\{e_1,\dots, e_{t_3}\}$.
Also, we may assume that $e_{t_3}\leq q-1$. 

Applying \rl{key} to the matrix $Q$, together with $e_{t_3}\leq q-1$ we obtain 
\begin{equation}\label{e:1}
g+f\leq \tdet Q+(q-1)t_4,
\end{equation}
where $f=f_1+\cdots+f_{t_3-1}+e_{t_3}$ and $g=g_1+\cdots+g_{t_3-1}+e_{t_3}+g_{t_3+1}+\cdots+g_{t_4}$.

Next note that for every $1\leq j\leq t_1$ we have $c_j\leq a_j$, by maximality of $T_Y$.
We also know that $i_j\leq q$ as it lies in the block $W$.   
Assume that we simultaneously have $c_j=a_j$ and $i_j=q$. Then if we  swap columns containing $c_j$ and $a_j$ we do not change $|T_Y|$ (as $c_j$ replaces $a_j$) but make $|T_Q|$ bigger (since $i_2=q$ replaces $e_{t_3}\leq q-1$), which contradicts our choice of $T_A$.  
Therefore, $c_j+i_j\leq a_j+q-1$ and, summing these up over $1\leq j\leq t_1$ we get
\begin{equation}\label{e:2}
c+i\leq \tdet Y+(q-1)t_1,
\end{equation}
where $c=c_1+\cdots+c_{t_1}$ and $i=i_1+\cdots+i_{t_1}$. Similarly, we have 
\begin{equation}\label{e:3}
d+h\leq \tdet Z+ (q-1)t_2,
\end{equation}
where $d=d_1+\cdots+d_{t_2}$ and $h=h_1+\cdots+h_{t_2}$.

Summing up \re{1}--\re{3} and using $c+f+h=mk$, $d+g+i=ml$, and \re{0} we get
$$mk+ml\leq \tdet Y+\tdet Z+\tdet Q+(q-1)(t_1+t_2+t_4)\leq \tdet A+(q-1)nl.
$$
Finally, this implies 
 \begin{eqnarray*}
 \tdet A\geq mk+ml-(q-1)nl&=&qnk+rk+qnl+rl-qnl+nl\\
 &=&qnk+r(k+l)+nl\geq qnk+nl\\
 &\geq& qnk+nk=nk(q+1).	
 \end{eqnarray*}
\end{pf}

Here is our main lower bound on the tropical determinant. In the next section we show that it is sharp.

\begin{Th}\label{T:main} Let $m=qn+r$ for $0\leq r< n$, and
$x,y$ as in \rd{xy}. Then 
$$L^{k,l}(m,n)\geq \min(nk(q+1), nkq+x+y).$$
\end{Th}

\begin{pf} As before, we can assume that
\[
A=\left(\begin{array}{cc}   
X&Y\\
Z&W
\end{array}\right)
\]
where $X$ is of size $t_1\times t_2$ and each entry of $W$ is at most $q$. If $t_1+t_2\geq nk$ then sum of dimensions of $W$ is
$$nk-t_1+nl-t_2\leq nl,
$$
so by \rc{marriage} there is a transversal in $A$ whose entries are at least $q+1$. Therefore, 
$$\tdet A\geq nk(q+1).$$

Now assume that $t_1+t_2< nk$. 
By \rc{Marriage} there exist  transversals  in  $Y$ and $Z$ whose entries are at least $q+1$. 
Thus, we can write $$\tdet Y\geq t_1(q+1)\quad\text{ and }\quad \tdet Z\geq t_2(q+1).$$
If we also have  $x+y\leq t_1+t_2$, then 
$$\tdet Y+\tdet Z+(nk-t_1-t_2)q\geq nkq+t_1+t_2\geq nkq+x+y.$$
The statement now follows from \rl{main}.

It remains to consider the case where   $t_1+t_2<nk$ and $t_1+t_2<x+y$.  
If we had $t_1\geq rk$ and $t_2\geq rl$, then \rl{lknr} and the definition of $x$ and $y$ 
would imply that $x+y\leq t_1+t_2$, which is not the case now.

If $t_1\leq rk$ and $t_2\leq rl$, then $(rk,rl)$ also satisfies the inequality in \rl{lknr}, so by \rd{xy}
we must have $x=rk$ and $y=rl$. On the other hand, since $t_1+t_2<nk$, by \rc{marriage}, 
every maximal transversal in $A$
contains an entry not exceeding $q$. Pick a maximal transversal and let $e$ be an entry in that transversal such that $e\leq q$.
By \rl{key}  we have
$$lm+km\leq \tdet A+lne\leq \tdet A+lnq,
$$
which implies 
$$\tdet A\geq nkq+kr+lr=nkq+x+y.
$$

Finally, assume that $t_1\geq rk$ and $t_2\leq rl$. As before, this implies that
$(t_1,rl)$ aslo  satisfies the inequality in \rl{lknr}, so by \rd{xy}
we must have  $x+y\leq t_1+rl$.
On the other hand, the row sums in $Z$ are bounded below by
$$ml-q(nl-t_2)=rl+qt_2,
$$
so by \rl{easy}, $\tdet Z\geq rl+qt_2$.
We have 
\begin{eqnarray*}
\tdet Y+\tdet Z+(nk-t_1-t_2)q&\geq& t_1(q+1)+rl+qt_2+(nk-t_1-t_2)q\\
&=&nkq+t_1+rl\geq nkq+x+y,
\end{eqnarray*}
and we are done by \rl{main}. The case $t_1\leq rk$ and $t_2\geq rl$ is similar.
\end{pf}

\section{Constructions}

To show that the bound in \rt{main} is sharp, we provide two constructions.

\begin{Prop}\label{P:constr1} There exists $A\in\mathcal D^{k,l}(m,n)$ such that $\tdet A\leq nk(q+1)$.
\end{Prop}
\begin{pf}
We describe how to construct such matrix $A\in\mathcal D^{k,l}(m,n)$ whose entries equal $q$ or $q+1$. Each row  of $A$ has $q+1$ repeated $rl$ times and each column  has $q+1$ repeated $rk$ times. To achieve this, in the  first row, place $q+1$ in the first $rl$ positions and fill in the remaining slots with $q$'s.  Let then each next row be a circular shift of the previous row by $rl$ slots.
In the resulting matrix we will have $rlnk$ entries equal $(q+1)$, so each column will contain  $rk$ $q+1$'s since we distributed them evenly among the columns. All the entries of $A$ are less than or equal to $q+1$, so $\tdet A\leq nk(q+1)$.
 Here is an example  of this construction with $m=7, n=5, k=1, l=3, q=1, r=2$.
 \[
A=\left(\begin{array}{ccccccccccccccc}   
2&2&2&2&2&2&1 &1 &1&1&1&1&1&1&1 \\
1&1&1&1&1&1&2 &2 &2&2&2&2&1&1&1 \\
2&2&2&1&1&1&1 &1 &1&1&1&1&2&2&2 \\
1&1&1&2&2&2&2 &2 &2&1&1&1&1&1&1 \\
1&1&1&1&1&1&1 &1 &1&2&2&2&2&2&2 
\end{array}\right).
\]
\end{pf}

Recall that $x$ and $y$ are described in \rd{xy}. Given that $x+y\leq nk$ we next explain how to construct a matrix with tropical determinant at most $nkq+x+y$.

\begin{Prop}\label{P:constr2} Let $x+y\leq nk$. Then there exists $A\in\mathcal D^{k,l}(m,n)$ such that $\tdet A\leq nkq+x+y$.
\end{Prop}
\begin{pf}
Let $A$ consist of four blocks $X$, $Y$, $Z$, and $W,$  that is,
\[
A=\left(\begin{array}{cc}   
X&Y\\
Z&W
\end{array}\right),
\]
where $X$ is of size $x\times y$. Let $W$ have all of its entries equal to $q$.  The only entries of $Z$ and $Y$ are $q$'s and $(q+1)$'s. 
We place $rk$ $(q+1)$'s in each column of $Y$ in a pattern similar to that of previous proposition. In the first column of $Y$ we place a string
of $rk$ $(q+1)$'s starting at first position (we can do this since $x\geq rk$) and fill in the remaining slots with $q$'s. In the second column we shift down this
string by $rk$ positions, circling around, if necessary. We repeat this  in every column of $Y$ starting with $(q+1)$'s in the position right after the one where we finished in the previous column.
We have distributed $\Sigma_Y=(xq+rk)(nl-y)$ (the sum of all the entries in~$Y$) evenly among the columns and as evenly as possible among the rows of $Y$. Let  $\Sigma_Y=ax+b$ be the result of dividing with remainder of 
$\Sigma_Y$ by~$x$, the number of rows in $Y$. Then the first $b$ rows of $Y$ have $q+1$ in  $a+1$ positions and the remaining $x-b$ rows have $q+1$ in $a$ positions.
Block  $Z$ is constructed in  a similar way, but we work with rows instead of columns. 

It remains to fill in block $X$. Its sum of entries is  $\Sigma_X=qxy+r(xl+yk)-klnr\geq 0$. We will distribute this sum as evenly as possible among the rows and columns of $X$. For this, we divide $\Sigma_X$ by $x$ with a remainder:
$\Sigma_X=cx+d$.  We want the bottom $d$ row sums of $X$ to  be equal to $c+1$ and the remaining row sums to be equal to $c$. 
For this we divide $c+1$  by~$y$ with remainder $c+1=ey +f$ and fill the last $f$ slots in the  bottom row 
of $X$ with $e+1$'s and make the remaining slots in the bottom row of $X$ equal to~$e$. Next row upward is a circular leftward shift of this row by~$f$. We continue with these circular shifts until we fill in the bottom $d$ rows of $X$. 
We fill the remaining rows of $X$ in a similar fashion: divide $c$ by $y$ with remainder $c=gy+h$, fill the
$h$ slots  in row $x-d$ with $(g+1)$'s (starting from where we stopped in the row below and going left) and so on.
% and continue in the same fashion, starting  with $h$ $(g+1)$'s in row $x-d$ to the left of the leftmost $e+1$ in the row below (or in the rightmost position if $f=0$ and there are no $(e+1)$'s). 
 In the resulting block $X$ the sum of entries $\Sigma_X$ is distributed as evenly as possible between rows and columns of~$X$. Moreover, the bigger row (resp. column) sums at the bottom (resp. rightmost) part of $X$. 
 We have also evenly distributed row sums in $Y$ and column sums in $Z$, so that bigger row sums in $Y$ are in the first rows of $Y$ and bigger column sums in $Z$ are in the first columns of $Z$. Hence first $x$ row sums and $y$ column sums of $A$ are equal to $ml$ and $mk$, respectively. 

Note that $e=g$ unless $c+1=ey$ and $f=0$, so the entries of $X$ differ from each other by at most 1.
Hence they are equal to $\floor*{\Sigma_X/xy}$  or $\floor*{\Sigma_X/xy}+1$, where the latter occurs only if $xy$ does not divide $\Sigma_{X}$ evenly. 
We have
$$\frac{\Sigma_X}{xy}=q+\frac{r(xl+yk)-lknr}{xy}\leq q
$$
since $xl+yk\leq lkn$ as
$$xl+yk\leq l(x+y)\leq lkn.
$$
This implies that each entry in $X$ is at most $q$. Hence for a maximal transversal of $A$ we can pick at most $x$ $(q+1)$'s in $Y$ and at most $y$ $(q+1)$'s in $Z$, so 
$$\tdet A\leq x(q+1)+y(q+1)+(nk-x-y)q=nkq+x+y,
$$
which completes the proof.
\end{pf}

Note that in the above argument we only used the conditions $x+y\leq nk$, $x\geq rk$, $y\geq rl$, and   
$qxy+r(xl+yk)\geq klnr$, but not the fact that $x+y$ is the smallest possible.
We next give an example of the above construction where this last assumption is dropped. This will allow us to  have the sum of entries in $X$ not too small, so that the construction of block $X$ can be better illustrated.

\begin{Ex}  Let $m=6, n=5,  k=2, l=3, q=1, r=1, x=5, y=4.$ Then each column of $Y$ has two 2's and three 1's and the sum of entries in $Y$ is $\Sigma_Y=(xq+rk)(ml-y)=77=5\cdot 15+2$, so we have row sums in first two rows equal to 16 and in the remaining three rows equal to 15. In block $Z$  the overall sum of entries equals
$$\Sigma_Z=(yq+rl)(mk-x)=35,$$ so first three  columns in $Z$ have columns sums equal to $9$, and the last column sum is 8. Next, 
$$\Sigma_X=qxy+r(xl+yk)-klnr=13,$$ so last  three row sums are 3, and first two are 2.  Therefore
 \[
A=\left(\begin{array}{cccc|ccccccccccc}   
1&0&0&1&2&1&2&1&1&2  &1&2&1&1&2 \\
0&1&1&0&2&1&1&2&1&2  &1&1&2&1&2 \\
1&1&0&1&1&2&1&2&1&1  &2&1&2&1&1 \\
1&0&1&1&1&2&1&1&2&1  &2&1&1&2&1 \\
0&1&1&1&1&1&2&1&2&1  &1&2&1&2&1 \\
\hline
2&2&2&1&1&1&1&1&1&1  &1&1&1&1&1 \\
2&2&1&2&1&1&1&1&1&1  &1&1&1&1&1 \\
2&1&2&2&1&1&1&1&1&1  &1&1&1&1&1 \\
1&2&2&2&1&1&1&1&1&1  &1&1&1&1&1 \\
2&2&2&1&1&1&1&1&1&1  &1&1&1&1&1 
\end{array}\right)
\]
is a matrix with tropical determinant at most $nkq+x+y=19$.
\end{Ex}

\rt{main} together with \rp{constr1} and \rp{constr2} imply  our main result.

\begin{Th}\label{T:mainshort}  Let $m=qn+r$ for $0\leq r< n$. Then 
$$L^{k,l}(m,n)=nkq+\min(nk, x+y),$$
 where $x,y$ are integers satisfying 
$$qxy+r(xl+yk)\geq klnr,\quad x\geq rk,\quad y\geq rl.$$ 

\end{Th}
\begin{Ex}\label{ex:3}
Let $m=6$, $n=5$, $k=1$, $l=2$, so $q=r=1$. It is easily seen that $x=y=2$ satisfy the inequalities
 $$xy+2x+y\geq 10,\quad x\geq 1,\quad y\geq 2,$$
 and have the smallest sum $x+y=4$ (see \rf{asym}). 
 \begin{figure}[h]\label{F:asym}
\includegraphics[scale=1.7]{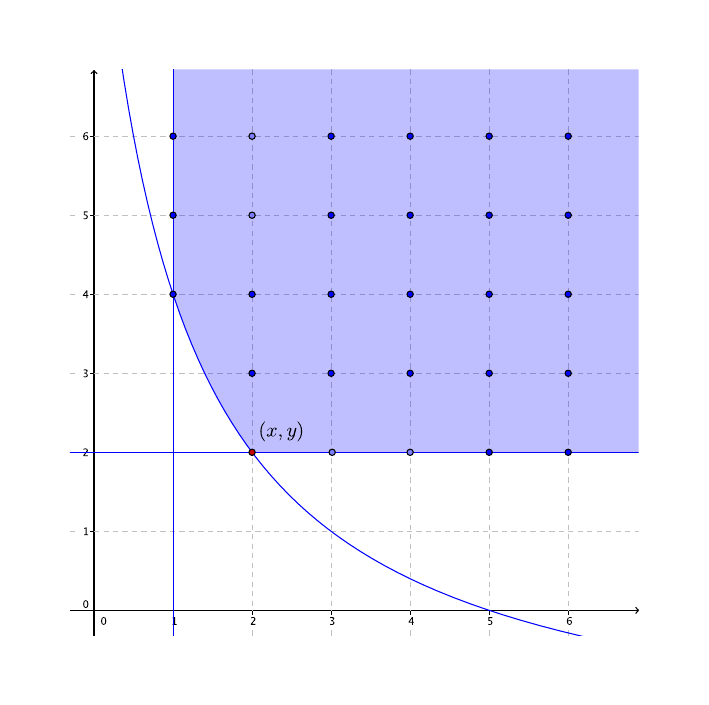}
\caption{Region in \rex{3} as defined in \rd{xy}.}
\end{figure}
 Also $x+y<nk=5$,  
so $L^{k,l}(m,n)=  nkq+x+y=9$, which is attained at the matrix in \rex{tdet9}.
 \end{Ex}

\section{Corollaries}
\begin{Cor}\label{C:0rem} If $r=0$, we have $L^{k,l}(m,n)= nkq$ and this value is attained at the matrix all of whose entries are equal to $q$.
\end{Cor}
\begin{pf}
Our conditions on $x$ and $y$ simplify to $qxy\geq 0$, $x\geq 0$, $y\geq 0$, where $x+y$ is the smallest possible, so $x=y=0$ and the main theorem implies $L^{k,l}(m,n)= nkq$.
Also, the matrix all of whose entries are equal to $q$ is in $\mathcal D^{k,l}(m,n)$ and its tropical determinant equals $nkq$.
\end{pf}

\begin{Cor} If $r\geq \frac{n}{q+2}$, we have $$L^{k,l}(m,n)=nkq+ \min(nk, r(k+l)).$$
\end{Cor}
\begin{pf}
The condition $r\geq \frac{n}{q+2}$ implies that $(rk, rl)$ satisfies the inequality \re{xy}. Therefore
$x=rk$ and $y=rl$, and the result follows from the main theorem.
\end{pf}

We now reformulate this corollary. 

\begin{Cor} 
If $r\geq\frac{n}{q+2}$ and $r\leq\frac{nk}{k+l}$ then $L^{k,l}(m,n)= nkq+r(k+l)$. 
If  $r\geq\frac{nk}{k+l}$ then $L^{k,l}(m,n)=nk(q+1)$. \\
\end{Cor}
\begin{pf}
To prove the second statement, we notice $x+y\geq rk+rl\geq nk$.
\end{pf}
 
 When $k=l=1$ the above corollary provides Theorems 2.3 and 2.4 from \cite{DHS}.
 Next we show how to recover the result of Theorem 3.3 from \cite{DHS}, which deals with
 the case $k=l=1$ and $r<\frac{n}{q+2}$.
  
% It remains to explore the case when $r<\frac{n}{q+2}$ and $r<\frac{nk}{k+l}$. 
% It turns out that under these assumptions there is no definite answer as to which one of $nk$ and $x+y$
% is smaller. For example, if $k=1, l=2, n=6, r=1,$ and  $q=1$, then $x+y=5$ and it is smaller than $nk=6$. On the other hand, if we change $l$ from 2 to 4, we get $x+y=7$, which is bigger than 
% $nk=6$.

%In the case when $k=l=1$ and $r<\frac{n}{q+2}$ (which implies $r<n/2$), we always get (as it was shown in \cite{DHS}) that $x+y\leq nk=n$. We now show how this result follows from our main theorem.
%

First, by definition, $x\geq r$, $y\geq r$ satisfy
$$qxy+r(x+y)-nr\geq 0,
$$
and $x+y$ is smallest possible. The region described by the above inequalities
is convex and symmetric with respect to the line $y=x$. It follows that
the region contains the segment joining points with coordinates $(r,n-r)$ and $(n-r,r)$ (see \rf{sym2}).
 \begin{figure}[h]\label{F:sym2}
\includegraphics[scale=2]{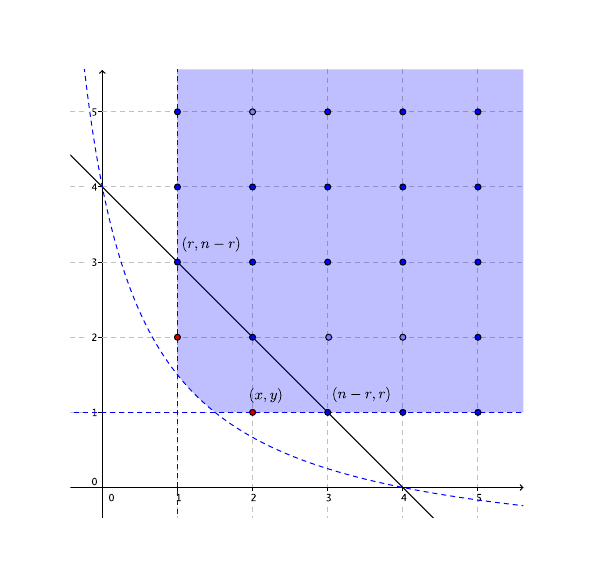}
\caption{Region described in \rd{xy} for $k=l=1$.}
\end{figure}
Therefore, any optimal solution $(x,y)$ satisfies $x+y\leq n$.
Furthermore, the minimum of $x+y$ is attained either when $y=x$ or $y=x+1$. Applying
\rt{mainshort} we get the statement of  \cite[Theorem 3.3]{DHS}.

\begin{Cor} Let $r<\frac{n}{q+2}$ and $x$ be the smallest positive integer 
satisfying at least one of the inequalities
\begin{enumerate}
\item $qx^2+2rx-nr\geq 0$,
\item $qx^2+(2r+q)x+r-nr\geq 0$.
\end{enumerate}
Then if $x$ satisfies (1) (and hence (2)), we have $L^{1,1}(m,n)=nq+2x$.
If $x$ satisfies (1) only, we have $L^{1,1}(m,n)=nq+2x+1$.
\end{Cor}

\section{Upper Bound on the  Tropical Determinant}
In this section we consider a version of the problem solved above where the maximum in the definition of the tropical determinant is replaced with  the minimum and we are interested in the
sharp upper bound of this tropical determinant on the transportation polytope.
Recall that $A=(a_{ij})$ is an $nk\times nl$ matrix where $\gcd(k,l)=1$, $k\leq l$, and $a_{ij}$ are nonnegative integers. The row sums of $A$ are equal to $ml$ and the column sums are equal to $mk$.
The set of all such matrices is denoted by $\cD^{k,l}(m,n)$.
As before, we divide $m$ by $n$ with remainder,  so $m=qn+r$, for $0\leq r<n$.

\begin{Def}  Let $A=(a_{ij})$ be an $s\times t$ matrix with $s\leq t$ and let $\mathcal{T}(A)$ be the set of its transversals. Define
the {\it tropical determinant} of a matrix $A=(a_{ij})$ to be
$$\tropdet(A)=\min_{T\in\mathcal{T}(A)}|T|.
$$
Denote its sharp upper bound over the set $\cD^{k,l}(m,n)$ by $U^{k,l}(m,n)$,
$$U^{k,l}(m,n)=\max_{A\in \cD^{k,l}(m,n)} \tropdet(A).
$$
\end{Def}

\begin{Th} $U^{k,l}(m,n)\leq \max(nkq,nkq+r(k+l)-nl)$.
\end{Th}
\begin{pf} Let $A\in\cD^{k,l}(m,n)$. Rearrange rows and columns  of $A$ so that the tropical determinant is equal to the sum of entries on the main diagonal of $A$ and the entries are non-decreasing along the main diagonal. That is,  $A$ is of the form
\begin{equation}
 A= \left( \begin{array}{cccccccc}
a_{11} &   & &      & a_{1t}  &&&\\
   & \ddots& &     & \vdots&&&\\
     &  &  a_{ii}&   & a_{it}&&&\\
   &   & & \ddots &\vdots &&&\\
a_{t1} & \dots &a_{ti} & \dots &   a_{tt}&a_{t\,t+1}&\dots&a_{ts}\\
\end{array} \right)
\label{e:matr}
\end{equation}
where $t=nk, s=nl$  $a_{11}\leq a_{22}\leq\cdots\leq a_{tt}$ and $\tropdet A=a_{11}+\cdots+a_{tt}$. Let us first suppose that $a_{tt}\leq q$. Then
$$\tropdet(A)=a_{11}+\cdots+a_{tt}\leq t\cdot a_{tt}=nkq.
$$
Next, let $a_{tt}\geq q+1$. Observe that 
$$a_{ti}+a_{it}\geq a_{ii}+a_{tt} \ \  {\rm for}\ \ i=1\dots t.
$$
since otherwise we could pick a smaller transversal. Also, for the same reason, $a_{tt}\leq a_{ti}$ for $i=t+1,\dots, s$. Adding  up all these inequalities over $i$ we get
$$a_{t1}+\cdots+a_{t\,t}+a_{t\, t+1}+a_{ts}+a_{1t}+\cdots+a_{t\, t}\geq \tropdet A +sa_{tt},
$$
and hence
$$mk+ml\geq \tropdet A+sa_{tt}\geq \tropdet A+nl(q+1),
$$
so
$$\tropdet A\leq m(k+l)-nl(q+1)=nkq+r(k+l)-nl.
$$
\end{pf}

\begin{Th} $U^{k,l}(m,n)= \max(nkq,nkq+r(k+l)-nl)$.
\end{Th} 
\begin{pf} Now it remains to construct matrices that reach the bound of the previous theorem. That is, for $r\leq nl/(k+l)$ we need to construct $A\in\cD^{k,l}(m,n)$ such that $\tropdet A\geq nkq$ and for
 $r\geq nl/(k+l)$ we need to construct $A\in\cD^{k,l}(m,n)$ such that $\tropdet A\geq nkq+r(k+l)-nl$.
The first task is easy. The entries of $A$ equal $q$ or $q+1$ with $rl$ $q+1$'s in each row, that are evenly distributed among the columns. That is, the first row of $A$ starts with $rl$ $q+1$'s and each next row is a circular 
shift by $rl$ of the previous row:
\[
A=\left(\begin{array}{cccccccccc}   
q+1&\dots&q+1&q&\dots&q&q&\dots &q&q \\
q&\dots&q&q+1&\dots&q+1&q&\dots &q& q\\
q&\dots&q&q&\dots&q&q+1&\dots &q+1& q\\
\dots&q+1&q&q&\dots&q&q&\dots &q& q+1\\
\dots&& & &\dots& & &\dots & & \\
\end{array}\right).
\]
There are $rlnk$ $q+1$'s in this matrix and since they are evenly distributed among the columns, each column contains $rlnk/nl=rk$ of them, so $A\in\cD^{k,l}(m,n)$. Since all the entries of $A$ are greater than or equal to $q$ we have
$\tropdet A\geq nkq$.

Let us next suppose that $r\geq nl/(k+l)$. Let  $A$ consist of four blocks
\[
A=\left(\begin{array}{cc}   
X&Y\\
Z&W
\end{array}\right),
\]
where $X$ is an $rk\times rl$ matrix of $q+1$'s, and the entries in $Y$ and $Z$ are all $q$. We fill in the remaining submatrix $W$ so that $A\in\cD^{k,l}(m,n)$. For this, we 
first make all entries of $W$ equal $q$.
We need to bring up the row sums in $W$ by $rl$ and the column sums by $rk$. For this, we divide $rl$ by $nl-rl$ with  remainder  to get $rl=(nl-rl)q'+r'$.
We increase the first $r'$ entries in the first row of $W$ by $q'+1$, and the remaining entries in this row by $q'$. 
 The second row of $W$ is a circular shift by $r'$ of the first row, and so on. Since we distributed $(ml-rlq)(nk-rk)$ as evenly as possible among the columns, the column sums in $W$ are
$$
\frac{(ml-rlq)(nk-rk)}{nl-rl}=mk-krq,
$$
so $A\in \cD^{k,l}(m,n)$. Note that all the entries in $W$ are greater than or equal to $q$. 

We have $nk-rk\leq rl$ and $nl-rl\leq rk$, so 
for a minimal transversal of $A$ we would need to pick $nk-rk$ entries
from $Z$, $nl-rl$ entries from~$Y$, and the remaining $r(k+l)-nl$ entries from $X$.
Therefore, $\tropdet A=nkq+r(k+l)-nl$.
  \end{pf}

\section{Acknowledgements.} Ivan Soprunov is partially supported by NSA Grant H98230-13-1-0279. Parts of this paper constituted Sailaja Gajula's Master's exit project at Cleveland State University written under the direction of Ivan Soprunov.

\end{document}